\documentclass[12pt]{article}
\usepackage{amsmath,amsfonts,amsthm,amssymb,eucal}
\DeclareMathOperator{\supp}{supp}

\begin{document}
\newtheorem{lem}{Lemma}
\newtheorem{teo}{Theorem}
\newtheorem{prop}{Proposition}
\pagestyle{plain}
\title{A Non-Archimedean Wave Equation}
\author{Anatoly N. Kochubei\footnote{Partially supported by
DFG under Grant 436 UKR 113/87/01, and by the Ukrainian Foundation 
for Fundamental Research, Grant 14.1/003.}\\ 
\footnotesize Institute of Mathematics,\\ 
\footnotesize National Academy of Sciences of Ukraine,\\ 
\footnotesize Tereshchenkivska 3, Kiev, 01601 Ukraine
\\ \footnotesize E-mail: \ kochubei@i.com.ua}
\date{}
\maketitle

\bigskip
\begin{abstract}
Let $K$ be a non-Archimedean local field with the normalized absolute value $|\cdot |$. 
It is shown that a ``plane wave'' $f(t+\omega_1 x_1+\cdots +\omega_nx_n)$, where $f$ is a
Bruhat-Schwartz complex-valued test function on $K$, $(t,x_1,\ldots ,x_n)\in K^{n+1}$, 
$\max\limits_{1\le j\le n}|\omega_j|=1$, satisfies, for any $f$, a certain homogeneous 
pseudo-differential equation, an analog of the classical wave equation. 
A theory of the Cauchy problem for this equation is developed.
\end{abstract}

\bigskip
{\bf Key words: }\ local field, plane wave, pseudo-differential equation, Cauchy problem

{\bf MSC 2000}. Primary: 11S80; 35S10. Secondary: 35L99. 

\bigskip
\section{INTRODUCTION}

Pseudo-differential equations for complex-valued functions defined 
on non-Archimedean local fields, in particular the field $\mathbb 
Q_p$ of $p$-adic numbers, are becoming increasingly important, 
both in view of rich mathematical structures involved in their 
studies and due to emerging applications; see the surveys in 
\cite{AKS1,K,KSA,Koz,Var,VVZ,Z}.

In most cases, pseudo-differential equations over $\mathbb 
Q_p^n$ with the symbols $|P(\xi_1,\ldots ,\xi_n)|_p^\alpha$, 
$\alpha >0$, where $P$ is a polynomial, were studied. The class of 
elliptic operators correspond to such polynomials $P$ that 
$P(\xi_1,\ldots ,\xi_n)\ne 0$ for $(\xi_1,\ldots ,\xi_n)\ne 0$. An 
equation of the Schr\"odinger type is obtained if there is a 
distinguished variable, say $\xi_1$, and $P(\xi_1,\ldots 
,\xi_n)=\xi_1-r(\xi_2,\ldots ,\xi_n)$ where $r$ is a $p$-adic 
quadratic form. It has been understood also that analogs of 
parabolic equations are evolution equations with a real time 
variable and $p$-adic spatial variables (this is connected with 
probabilistic applications; see \cite{K} and references therein). 
It seemed natural to interpret the case where $P(\xi_1,\ldots 
,\xi_n)=\xi_1^2-r(\xi_2,\ldots ,\xi_n)$ with an anisotropic 
quadratic form $r$, as the case of hyperbolic equations. However 
the results available for this case (see \cite{K}) are quite 
scant. In particular, no well-posed problems for such equations 
have been identified.

In this paper we propose an alternative approach. Instead of a 
formal resemblance in the definition of an equation, we proceed 
from its properties. Let us call a function $u(t,x_1,\ldots 
,x_n):\ \mathbb Q_p^{n+1}\to \mathbb C$ {\it a plane wave}, if 
\begin{equation}
u(t,x_1,\ldots ,x_n)=f(t+\omega_1 x_1+\cdots +\omega_nx_n)
\end{equation}
where $f$ belongs to the Bruhat-Schwartz space $\mathcal D(\mathbb 
Q_p)$ of test functions, $\omega_1,\ldots ,\omega_n\in \mathbb 
Q_p$, $\max\limits_{1\le j\le n}|\omega_j|_p=1$ (in fact, we will 
consider not only $\mathbb Q_p$ but arbitrary non-Archimedean 
local fields; see below).

We will show that every function (1) is a solution of the equation
\begin{equation}
D_t^\alpha u-D_x^{\alpha ,n}u=0
\end{equation}
where $D^\alpha$ is Vladimirov's fractional differentiation 
operator, that is a pseudo-differential operator with the symbol 
$|\xi |_p^\alpha$, while $D^{\alpha ,n}$ is a pseudo-differential 
operator of $n$ variables with the symbol $\max\limits_{1\le j\le 
n}|\xi_j|_p^\alpha$, $\alpha >0$ is an arbitrary number.

The equation (2) with $n=1$ was mentioned in \cite{V03} as an 
example of the following pathology. Consider the equation for a 
related fundamental solution $E$:
\begin{equation}
D_t^\alpha E-D_x^{\alpha ,n}E=\delta
\end{equation}
where $E$ belongs to some class of distributions, on which the 
operators are defined, with the usual relations between them and 
the Fourier transform. Then, performing the Fourier transform we 
obtain the contradictory identity $\left( |\tau |_p^\alpha -
\max\limits_{1\le j\le n}|\xi_j|_p^\alpha \right) 
\widetilde{E}(\tau ,\xi_1,\ldots ,\xi_n)=1$ where the left-hand 
side vanishes on the open set
$$
\left\{ 0\ne (\tau ,\xi_1,\ldots ,\xi_n)\in \mathbb Q_p^{n+1}:\
|\tau |_p=\max\limits_{1\le j\le n}|\xi_j|_p\right\} .
$$

Therefore the fundamental solution cannot exist, and one cannot 
expect any reasonable behavior of an inhomogeneous equation 
associated with (2). On the other hand, the set of solutions of 
the one-dimensional equation $D_t^\alpha u=\lambda u$ is 
infinite-dimensional (see \cite{K,VVZ}). Thus, at the first sight, 
the equation (2) does not look as an evolution equation with the 
``time'' variable $t$.

Nevertheless, in this paper we prove the existence and uniqueness 
of solutions for some analogs of the Cauchy problem for the 
equation (2) in the class of {\it radial} functions, that is those 
depending (in the variable $t$) only on $|t|_p$. On this class, 
the operator $D^\alpha$ becomes a counterpart of the 
Caputo-Dzhrbashyan regularized fractional derivative appearing in 
fractional evolution equations of real analysis (see \cite{EIK}). 
Moreover, the above connection with plane waves, together with the 
inversion formula for the Radon transform, available in the 
non-Archimedean case too \cite{Ch1}, leads to a formula for 
solutions and an analog of the Huygens principle.

\section{Preliminaries}

{\bf 2.1. Local fields.} Let $K$ be a non-Archimedean local field, 
that is a non-discrete totally disconnected locally compact 
topological field. It is well known that $K$ is isomorphic either 
to a finite extension of the field $\mathbb Q_p$ of $p$-adic 
numbers (if $K$ has characteristic 0), or to the field of formal 
Laurent series with coefficients from a finite field $\mathbb 
F_q$, if $K$ has characteristic $p\ne 0$; in this case $q=p^\nu$, 
$\nu \in \mathbb N$. For a summary of main notions and results 
regarding local fields see, for example, \cite{K}.

Any local field is endowed with an absolute value $|\cdot |$, 
such that $|x|=0$ if and only if $x=0$, $|xy|=|x|\cdot |y|$, 
$|x+y|\le \max (|x|,|y|)$. Denote $O=\{ x\in K:\ |x|\le 1\}$,
$P=\{ x\in K:\ |x|<1\}$, $U=O\setminus P$. $O$ is a subring of $K$ 
called the ring of integers, $P$ is an ideal in $O$ called the 
prime ideal; the multiplicative group $U$ is called the group of 
units. The ideal $P$ contains an element $\beta$, such that 
$P=\beta O$ (a prime element). The quotient ring $O/P$ is  
actually a finite field; denote by $q$ its cardinality. We will 
always assume that the absolute value is 
normalized, that is $|\beta |=q^{-1}$. The normalized absolute 
value $|\cdot |$ takes the values $q^N$, $N\in \mathbb Z$.

If $K=\mathbb Q_p$ ($p$ is a prime number), that is a completion 
of the field $\mathbb Q$ of rational numbers with respect to the 
absolute value
$$
|x|_p=p^{-\nu }\quad \text{for $x=p^\nu \frac{m}n$,}
$$
where $\nu ,m,n\in \mathbb Z$, and $m,n$ are prime to $p$, then 
$\beta =p$ ($p$ is seen as an element of $\mathbb Q_p$) and $q=p$ 
(as a natural number).

Returning to a general local field $K$, denote by $S\subset O$ a 
complete system of representatives of the residue classes from 
$O/P$. Then any nonzero element $x\in K$ admits the canonical 
representation in the form of the convergent series
$$
x=\beta^{-n}\left( x_0+x_1\beta +x_2\beta^2+\cdots \right)
$$
where $n\in \mathbb Z$, $|x|=q^n$, $x_j\in S$, $x_0\notin P$. For 
$K=\mathbb Q_p$, one may take $S=\{ 0,1,\ldots ,p-1\}$.

The additive group of any local field is self-dual, that is if 
$\chi$ is any non-constant complex-valued additive character of 
$K$, then any other additive character can be written as 
$\chi_a(x)=\chi (ax)$, $x\in K$, for some $a\in K$. See \cite{K} 
for an explicit description of the character $\chi$ used in 
harmonic analysis on local fields (``the canonical additive 
character''). In particular, it is assumed that $\chi$ is a rank 
zero character, that is $\chi (x)\equiv 1$ for $x\in O$, while 
there exists such an element $x_0\in K$ that $|x_0|=q$ and $\chi 
(x_0)\ne 1$.

The above duality is used in the definition of the Fourier 
transform over $K$. Denoting by $dx$ the Haar measure on the 
additive group of $K$ (normalized in such a way that the measure 
of $O$ equals 1) we write
$$
\widetilde{f}(\xi )=\int\limits_K\chi (x\xi )f(x)\,dx,\quad \xi 
\in K,
$$
where $f$ is a complex-valued function from $L_1(K)$. As usual, the Fourier 
transform $\mathcal F$ can be extended from $L_1(K)\cap L_2(K)$ to a 
unitary operator on $L_2(K)$. If $\mathcal F f=\widetilde{f}\in L_1(K)$, we 
have the inversion formula
$$
f(x)=\int\limits_K\chi (-x\xi )\widetilde{f}(\xi )\,d\xi .
$$

Similarly, if $f:\ K^n\to \mathbb C$, we write
$$
\widetilde{f}(\xi_1,\ldots ,\xi_n)=\int\limits_{K^n}\chi 
(x_1\xi_1+\cdots +x_n\xi_n)f(x_1,\ldots ,x_n)\,dx_1\ldots dx_n.
$$
The inversion formula is then
$$
f(x_1,\ldots ,x_n)=\int\limits_{K^n}\chi (-x_1\xi_1-\cdots -x_n\xi_n)
\widetilde{f}(\xi_1,\ldots ,\xi_n)\,d\xi_1\ldots d\xi_n.
$$

\medskip
{\bf 2.2. Spaces of test functions and distributions.} A function 
$f:\ K\to \mathbb C$ is called locally constant, if there exists 
such an integer $l$ that for any $x\in K$
$$
f(x+x')=f(x), \quad \text{if $|x'|\le q^{-l}$}.
$$
The smallest number $l$ with this property is called the exponent 
of local constancy of the function $f$.

Denote by $\mathcal D(K)$ the set of all locally constant 
functions with compact supports. $\mathcal D(K)$ is a vector space 
over $\mathbb C$. In order to furnish it with a topology, consider 
a subspace $\mathcal D_N^l\subset \mathcal D(K)$ of functions with 
supports in the ball
$$
B_N=\left\{ x\in K:\ |x|\le q^N\right\} ,\quad n\in \mathbb Z,
$$
and the exponents of local constancy $\le l$.

The space $\mathcal D_N^l$ is finite-dimensional; thus it has a 
natural topology induced from $\mathbb C$. Then we set 
$$
\mathcal D_N=\varinjlim_{l\to \infty }\mathcal D_N^l,
$$
and define the topology in $\mathcal D(K)$ as the inductive
limit topology, that is
$$
\mathcal D(K)=\varinjlim_{N\to \infty }\mathcal D_N.
$$
The strong conjugate space $\mathcal D'(K)$ is called the space of 
Bruhat-Schwartz distributions.

The operation of the Fourier transform preserves the space $\mathcal D(K)$
or $\mathcal D(K^n)$ (this property contrasts the Archimedean 
case). Therefore the Fourier transform of a distribution defined, 
by duality, just as for distributions from $\mathcal S'(\mathbb 
R^n)$, acts continuously on $\mathcal D'(K)$, resp. $\mathcal 
D'(K^n)$. As in the case of $\mathbb R^n$, there exists a 
well-developed theory of distributions over local fields including 
such topics as convolution, direct product, homogeneous 
distributions etc. Note in particular that a function $|x|^{\alpha 
-1}$, $\Re \ \alpha >0$, admits an analytic continuation in $\alpha$ 
to a meromorphic distribution
\begin{equation}
\left\langle |x|^{\alpha -1},\varphi \right\rangle =\int\limits_K
 |x|^{\alpha -1}[\varphi (x)-\varphi (0)]\,dx,\quad \varphi \in 
\mathcal D(K),
\end{equation}
$\Re \alpha >0$ (see Sect. VIII in \cite{VVZ} for $K=\mathbb Q_p$. 
The general case is completely similar). See \cite{AKS2,GGP,K,VVZ} 
for further details.

Below we will often use the subspaces of $\mathcal D(K^n)$,
$$
\Psi (K^n)=\left\{ \psi \in \mathcal D(K^n):\ \psi (0)=0\right\} ,
$$
$$
\Phi (K^n)=\left\{ \varphi \in \mathcal D(K^n):\ \int\limits_{K^n}
\varphi (x)\,d^nx=0\right\} ,
$$
introduced in \cite{AKS1}. The space $\Phi (K^n)$ is called the 
Lizorkin space of test functions of the second kind; it is a 
non-Archimedean counterpart of a space of test functions on 
$\mathbb R^n$ proposed by Lizorkin (\cite{L}; see also \cite{S}). 
The conjugate space $\Phi'(K^n)$ is called  the 
Lizorkin space of distributions of the second kind. The most 
important property of these spaces is that the Fourier transform 
$\mathcal F$ is a linear isomorphism from $\Psi (K^n)$ onto $\Phi (K^n)$, 
thus also from $\Phi'(K^n)$ onto $\Psi'(K^n)$. At the same time, 
$\mathcal F$ can be considered as a linear isomorphism from $\Phi (K^n)$ to
$\Psi (K^n)$.

\medskip
{\bf 2.3. Pseudo-differential operators.} The simplest and best 
studied pseudo-differential operator, acting on complex-valued 
functions over $K$, is the fractional differentiation operator 
$D^\alpha$, $\alpha >0$, whose deep investigation was initiated by 
Vladimirov (see \cite{VVZ}). It is defined as
$$
\left( D^\alpha \varphi \right) (x)=\mathcal F^{-1}\left[ |\xi |^\alpha 
(\mathcal F (\varphi ))(\xi )\right] (x),\quad \varphi \in \mathcal D(K).
$$
Note that $D^\alpha$ does not act on the space $\mathcal D(K)$, 
since the function $\xi \mapsto |\xi |^\alpha$ is not locally 
constant. We can assert, for example, that $D^\alpha :\ \mathcal 
D(K) \to L_2(K)$, and the closure of this operator is self-adjoint 
on $L_2(K)$. On the other hand, $D^\alpha :\ \Phi (K)\to \Phi (K)$ 
and  $D^\alpha :\ \Phi'(K)\to \Phi'(K)$; see \cite{AKS1}.

Similarly, for $x=(x_1,\ldots ,x_n)\in K^n$, set 
$\| x\|=\max\limits_{1\le j\le n}|x_j|$. The pseudo-differential operator 
$D^{\alpha ,n}:\ \mathcal D(K^n) \to L_2(K^n)$ is given by the 
expression
$$
\left( D^{\alpha ,n}\varphi \right) (x)=\mathcal F^{-1}\left[ \|\xi \|^\alpha 
(\mathcal F (\varphi ))(\xi )\right] (x),\quad \varphi \in \mathcal D(K^n).
$$
We have $D^{\alpha ,n}:\ \Phi (K^n)\to \Phi (K^n)$ 
and  $D^{\alpha ,n}:\ \Phi'(K^n)\to \Phi'(K^n)$.

An important property of these operators is the possibility to get 
rid of the Fourier transform and represent them as hyper-singular 
integral operators. For any $u\in \mathcal D(K)$,
\begin{equation}
\left( D^\alpha u\right) (x)=\frac{1-q^\alpha }{1-q^{-\alpha 
-1}}\int\limits_K |y|^{-\alpha -1}[u(x-y)-u(x)]\,dy;
\end{equation}
see \cite{K,VVZ}. The expression in the right-hand side of (5) 
makes sense for wider classes of functions, for example, for all 
bounded locally constant functions.

Similarly, if $u\in \mathcal D(K^n)$, then
\begin{equation}
\left( D^{\alpha ,n}u\right) (x)=\frac{1-q^\alpha }{1-q^{-\alpha 
-n}}\int\limits_{K^n}\|y\|^{-\alpha -n}[u(x-y)-u(x)]\,d^ny
\end{equation}
(see \cite{T68,AKS1,VVZ}).

\medskip
\begin{lem}
If $u$ is a bounded locally constant function on $K^n$, then the 
distribution $D^{\alpha ,n}u\in \Phi'(K^n)$ coincides with the 
function (6).
\end{lem}

\medskip
{\it Proof}. Let $\varphi \in \Phi (K^n)$. Then
$$
\left\langle D^{\alpha ,n}u,\varphi \right\rangle =\frac{1-q^\alpha }{1-q^{-\alpha 
-n}}\int\limits_{K^n}u(x)\,d^nx\int\limits_{K^n}\frac{\varphi 
(y)-\varphi (x)}{\| x-y\|^{\alpha +n}}\,d^ny.
$$
Let $\theta >0$ be so small that $u(x)=u(y)$ and $\varphi 
(x)=\varphi (y)$ if $\|x-y\| <\theta$. Denote
$$
C_\alpha =\frac{1-q^\alpha }{1-q^{-\alpha 
-n}}\int\limits_{\|y\|\ge \theta}\frac{dy}{\|y\|^{\alpha +n}}.
$$
Then
\begin{multline*}
\left\langle D^{\alpha ,n}u,\varphi \right\rangle =\frac{1-q^\alpha }{1-q^{-\alpha 
-n}}\int\limits_{K^n}u(x)\,d^nx\int\limits_{\|x-y\|\ge \theta}
\frac{\varphi (y)}{\| x-y\|^{\alpha +n}}\,d^ny-C_\alpha
\int\limits_{K^n}u(x)\varphi (x)\,d^nx\\
=\frac{1-q^\alpha }{1-q^{-\alpha -n}}\int\limits_{K^n}\varphi 
(y)\,d^ny\int\limits_{\|x-y\|\ge \theta}
\frac{u(x)-u(y)}{\| x-y\|^{\alpha +n}}\,d^nx=\langle \varphi ,\psi 
\rangle
\end{multline*}
where $\psi$ is the right-hand side of (6), as desired. $\qquad 
\blacksquare$

\medskip
{\bf 2.4. The Radon transform.} Let $\varphi \in \mathcal D(K^n)$, 
$n\ge 2$. The Radon transform $\check \varphi (\xi ,s)$, where 
$\xi \in K^n$, $\xi \ne 0$, $s\in K$, is defined by the relation
$$
\check \varphi (\xi ,s)=\int\limits_{\xi \cdot x=s}\varphi 
(x)\,d\omega_{\xi ,s}(x)
$$
(see \cite{Ch1}) where $\omega_{\xi ,s}$ is such a measure on the 
hyperplane $\xi \cdot x=s$ (we write $\xi \cdot x=\xi_1x_1+\cdots 
+\xi_nx_n$) that for any $\psi \in \mathcal D(K^n)$,
$$
\int\limits_Kds\int\limits_{\xi \cdot x=s}\psi 
(x)\,d\omega_{\xi ,s}(x)=\int\limits_{K^n}\psi (x)\,dx.
$$

The function $\check \varphi$ possesses the following properties. 
It is homogeneous of degree -1 in $\xi$ and $s$, that is $\check 
\varphi (\sigma \xi ,\sigma s)=|\sigma |^{-1}\check \varphi (\xi 
,s)$, for any $\sigma \in K\setminus \{0\}$. Next, $\check \varphi 
(\xi ,s)=0$, if the expression $|s|\cdot \| \xi \|^{-1}$ is sufficiently 
large. The function $\check \varphi$ is jointly locally constant 
in $\xi$ and $s$. Finally, the integral $\int\limits_K\check \varphi 
(\xi ,s)\,ds$ does not depend on $\xi$. Note that the above 
properties of a function of  $\xi$ and $s$ are also sufficient for 
such a function to be the Radon transform of some function $\varphi \in \mathcal 
D(K^n)$.

In order to find a connection between the Radon and Fourier 
transforms (similar to the well-known one for the case of $\mathbb 
R^n$ \cite{Hel}), we write
$$
\widetilde{\varphi }(s\xi )=\int\limits_{K^n}\varphi (x)\chi (s(x\cdot
\xi ))\,d^nx=\int\limits_Kdr\int\limits_{\xi \cdot x=r}\varphi 
(x)\chi (sr)\,d\omega_{\xi ,r}(x)
=\int\limits_K\chi (sr)\check \varphi (\xi ,r)\,dr,
$$
and it follows from the Fourier inversion formula that
\begin{equation}
\check \varphi (\xi ,r)=\int\limits_K\chi (-sr)\widetilde{\varphi }(s\xi 
)\,ds.
\end{equation}

The inversion formula for the non-Archimedean Radon transform is 
as follows \cite{Ch1}:
\begin{equation}
\varphi (x)=\frac{1-q^{n-1}}{(1-q^{-1})(1-q^{-n})}\int\limits_{\| 
\eta \|=1}\left\langle |s|^{-n},\check \varphi (\eta ,s+\eta \cdot 
x)\right\rangle \,d^n\eta ,\quad x\in K^n,
\end{equation}
where the distribution $|s|^{-n}$ is understood in the sense of 
(4). Substituting (4) into (8) and comparing with (5) we can write 
the inversion formula in the following form:
\begin{equation}
\varphi (x)=\frac{1}{1-q^{-1}}\int\limits_{\| \eta \|=1}\left. \left( 
D_s^{n-1}\check \varphi (\eta ,s+\eta \cdot 
x)\right) \right|_{s=0}d^n\eta .
\end{equation}

The identity (9) can be proved directly, by substituting (7) and 
calculating the integrals.

If $n=1$, we define the Radon transform by the formula (7). It is 
easy to check that the inversion formula (9) remains valid for 
this case too, in the form
$$
\varphi (x)=(1-q^{-1})^{-1}\int\limits_{|\eta |=1}\check \varphi (\eta 
,\eta x)\,d\eta .
$$

\section{Radial Eigenfunctions}

{\bf 3.1. $L_2$-solutions.} Let $u(x)=\psi (|x|)\in L_2(K)$,
\begin{equation}
D^\alpha u=\lambda u,\quad \lambda =q^{\alpha N},\ N\in \mathbb Z,
\end{equation}
and $u$ is not identically zero.

Let us apply the Fourier transform to both sides of (10). We get
\begin{equation}
\left( |\xi |^\alpha -q^{\alpha N}\right) \widetilde{u}(\xi )=0 
\quad \text{for all $\xi \in K$}.
\end{equation}
It follows from (11) that the inequality $\widetilde{u}(\xi )\ne 
0$ is possible only for $|\xi |=q^N$. Since $u$ is a radial 
function, $\widetilde{u}$ also possesses this property 
\cite{K,VVZ}. Therefore
\begin{equation}
\widetilde{u}(\xi )=\begin{cases}
c, & \text{if $|\xi |=q^N$;}\\
0, & \text{if $|\xi |\ne q^N$,}
\end{cases}\quad c\ne 0.
\end{equation}
By the Fourier inversion and the well-known integration formula 
(see \cite{K,VVZ}), we get
\begin{equation}
u(x)=\begin{cases}
cq^N(1-q^{-1}), & \text{if $|x|\le q^{-N}$;}\\
-cq^{N-1}, & \text{if $|x|=q^{-N+1}$;}\\
0, & \text{if $|x|>q^{-N+1}$.}\end{cases}
\end{equation}
It is easily seen from (12) or (13) that $u\in \Phi (K)$.

The only radial eigenfunction $u$ with $u(0)=1$ (an analog of the 
function $t\to e^{-\lambda t}$, $t\in \mathbb R$) corresponds to 
$c=q^{-N}(1-q^{-1})^{-1}$. On the other hand, if $u(0)=0$, then 
$c=0$.

\medskip
{\bf 3.2. Generalized solutions.} Let us consider solutions $u\in 
\Phi'(K)$ of the equation (10). It is natural to call a 
distribution $u\in \Phi'(K)$ radial (or spherically symmetric), if, 
for any $\omega \in K$, $|\omega |=1$, and any $\varphi \in \Phi 
(K)$
$$
\langle u,\varphi_\omega \rangle =\langle u,\varphi \rangle
$$
where $\varphi_\omega (x)=\varphi (\omega x)$. In a similar way, 
we define a radial distribution from $\Psi'(K)$. It is easy to 
check that the Fourier transform maps a radial distribution from 
$\Phi'(K)$ to a radial distribution from $\Psi'(K)$.

\medskip
\begin{prop}
If a radial distribution $u\in \Phi'(K)$ satisfies the equation 
(10), then it coincides, for some $c\in \mathbb C$, with the function 
(13).
\end{prop}

\medskip
{\it Proof.} By definition of a generalized solution, we have
$$
\left\langle u,D^\alpha \varphi \right\rangle =\lambda \langle u,\varphi \rangle
\quad \text{for any $\varphi \in \Phi (K)$}.
$$
Writing $\varphi =\mathcal F^{-1}\psi$, $\psi \in \Psi (K)$, we see that
$$
\left( D^\alpha \varphi \right) (x)=\mathcal F^{-1}_{\xi \to x}\left( |\xi |^\alpha
\psi (\xi )\right) .
$$

The function $\xi \to |\xi |^\alpha \psi (\xi )$ belongs to $\Psi (K)$. 
Therefore, considering $\mathcal F$ as an operator from $\Phi'(K)$ to 
$\Psi'(K)$, we may write
$$
\left\langle u,D^\alpha \varphi \right\rangle =\left\langle 
(\mathcal F u)(\xi ),|\xi |^\alpha \psi (\xi )\right\rangle = 
\left\langle |\xi |^\alpha (\mathcal F u)(\xi ),\psi (\xi )\right\rangle ,
$$
so that we come to the equality (11) where this time 
$\widetilde{u}=\mathcal F u\in \Psi'(K)$, and the multiplication by $|\xi 
|^\alpha -q^{\alpha N}$ is understood in the distribution sense. 
Thus, for any $\psi \in \Psi (K)$,
$$
\left\langle (\mathcal F u)(\xi ),\left( |\xi |^\alpha -q^{\alpha N}\right)
\psi (\xi )\right\rangle =0.
$$

On the sphere $\{ \xi \in K:\ |\xi |=q^l\}$, $l\ne N$, the set of 
functions $\xi \mapsto \left( |\xi |^\alpha -q^{\alpha N}\right)
\psi (\xi )$ runs the set of restrictions of all the functions 
from $\Psi (K)$. Therefore the restriction of the distribution 
$\mathcal F u$ to such a sphere equals zero, so that $\mathcal F u$ is concentrated on 
the sphere $S_N=\{ \xi \in K:\ |\xi |=q^N\}$.

The set of restrictions to $S_N$ of functions from $\Psi (K)$ 
coincides with $\mathcal D(S_N)=\varinjlim\limits_{l\to \infty}\mathcal 
D_l(S_N)$ where $\mathcal D_l(S_N)$ is the set of complex-valued 
functions on $S_N$ with the exponents of local constancy $\le l$. 
The space $\mathcal D_l(S_N)$ is finite-dimensional; its basis can 
be constructed from the functions 
$\delta_{\sigma_0,\sigma_1,\ldots ,\sigma_{N+l-1}}(t)$ 
($\sigma_0,\sigma_1,\ldots ,\sigma_{N+l-1}\in S$, $\sigma_0\notin 
P$), which equal 1 on elements $t\in S_N$ with the canonical 
representations $t=\beta^{-N}\left( \sigma_0+\sigma_1\beta +\cdots 
+\sigma_{N+l-1}\beta^{N+l-1}\right) +O(\beta^l)$, and 0 on all 
other $t\in S_N$.

Denote
$$
c_{\sigma_0,\sigma_1,\ldots ,\sigma_{N+l-1}}=\left\langle \mathcal F u,
\delta_{\sigma_0,\sigma_1,\ldots ,\sigma_{N+l-1}}\right\rangle .
$$
The ratio of any two elements $\beta^{-N}\left( \sigma_0+\sigma_1\beta +\cdots 
+\sigma_{N+l-1}\beta^{N+l-1}\right)$ belongs to the group of units 
$U$. The transformation of one of the functions $\delta_{\sigma_0,\sigma_1,\ldots 
,\sigma_{N+l+1}}$ into another (with the same $l$) is implemented 
by the multiplication of the argument by the appropriate ratio. 
Since $\mathcal F u$ is a radial distribution, we find that $c_{\sigma_0,\sigma_1,\ldots 
,\sigma_{N+l-1}}$ depends only on $l$, say
$$
c_{\sigma_0,\sigma_1,\ldots ,\sigma_{N+l-1}}=c'_{l-1}.
$$
At the same time,
$$
\sum\limits_{\sigma_{N+l}\in S}\delta_{\sigma_0,\sigma_1,\ldots ,
\sigma_{N+l-1},\sigma_{N+l}}=\delta_{\sigma_0,\sigma_1,\ldots 
,\sigma_{N+l-1}},
$$
whence $c'_{l-1}=qc'_l$ and $c_l'=c_0'q^{-l}$, $c'_0\in \mathbb 
C$. Thus, we have found that
\begin{equation}
\left\langle \mathcal F u,\delta_{\sigma_0,\sigma_1,\ldots 
,\sigma_{N+l}}\right\rangle =c_0'q^{-l}
\end{equation}
for all $l$.

Meanwhile, the integral
$$
\int\limits_{|t|=q^N}\delta_{\sigma_0,\sigma_1,\ldots 
,\sigma_{N+l}}(t)\,dt
$$
equals $q^{N-(N+l)-1}=q^{-l-1}$ (see Sect. IV in \cite{VVZ}). 
Together with (14), this shows that the restriction of the 
distribution $\mathcal F u$ to the sphere $S_N$ equals a constant; outside 
$S_N$, $\mathcal F u$ equals 0. Thus, $\mathcal F u$ has the form (12), so that $u$ 
coincides with the function (13).$\qquad \blacksquare$

\section{Plane Waves}

Following a classical pattern we call a function
\begin{equation}
F(t,x)=f(t+\omega_1x_1+\cdots +\omega_nx_n),
\end{equation}
$t\in K$, $(x_1,\ldots ,x_n)\in K^n$, where $\| (\omega_1,\ldots 
,\omega_n)\| =1$, $f\in \mathcal D(K)$, {\it a non-Archimedean 
plane wave}.

\medskip
\begin{prop}
For any $\alpha >0$, a non-Archimedean plane wave (15) satisfies 
the equation
\begin{equation}
D_t^\alpha F-D_x^{\alpha ,n}F=0.
\end{equation}
\end{prop}

\medskip
{\it Proof.} Suppose that $n\ge 2$ (in the case $n=1$ the validity 
of (16) is checked in a straightforward way). Let us compute $D_x^{\alpha 
,n}F$. By the definition of $D_x^{\alpha ,n}F$,
\begin{multline*}
\left( D_x^{\alpha ,n}F\right) (t,x)
=\frac{1-q^\alpha }{1-q^{-n-\alpha }}\int\limits_{K^n}\left( 
\max\limits_j|y_j|\right)^{-n-\alpha }\\
\times\left[ f\left( 
t+\sum\limits_{j=1}^n\omega_jx_j-\sum\limits_{j=1}^n\omega_jy_j\right)
-f\left( t+\sum\limits_{j=1}^n\omega_jx_j\right) 
\right]\,dy_1\ldots dy_n.
\end{multline*}

Since $\max\limits_j|\omega_j|=1$, we can choose an index $j_0$ in 
such a way that $\left| \omega_{j_0}\right| =1$. Suppose for 
simplicity that $|\omega_1|=1$. Let us perform the change of 
variables
$$
\eta_1=\sum\limits_{j=1}^n\omega_jy_j,\ \eta_2=y_2,\ldots 
,\eta_n=y_n.
$$
Obviously, $\max\limits_j|\eta_j|\le \max\limits_j|y_j|$. On the 
other hand,
$$
y_1=\frac{1}{\omega_1}\left( \eta_1-\omega_2\eta_2-\cdots 
-\omega_n\eta_n\right) ,
$$
whence $\max\limits_j|y_j|\le \max\limits_j|\eta_j|$, so that
$$
\max\limits_j|y_j|=\max\limits_j|\eta_j|.
$$

The Jacobian of the transformation $(y_1,\ldots ,y_n)\mapsto 
(\eta_1,\ldots ,\eta_n)$ equals
$$
\begin{vmatrix}
\omega_1 & \omega_2 & \dots & \omega_n\\
0 & 1 & \dots & 0\\
\hdotsfor{4} \\
0 & 0 & \dots & 1
\end{vmatrix}
$$
and belongs to $U$. We have
\begin{multline*}
\left( D_x^{\alpha ,n}F\right) (t,x)
=\frac{1-q^\alpha }{1-q^{-n-\alpha }}\int\limits_{K^n}\left( 
\max\limits_j|\eta_j|\right)^{-n-\alpha }\\
\times\left[ f\left( 
t+\sum\limits_{j=1}^n\omega_jx_j-\eta_1\right)
-f\left( t+\sum\limits_{j=1}^n\omega_jx_j\right) 
\right]\,d\eta_1\ldots d\eta_n\\
=\frac{1-q^\alpha }{1-q^{-n-\alpha }}\int\limits_K
\left[ f\left( 
t+\sum\limits_{j=1}^n\omega_jx_j-\eta_1\right)
-f\left( t+\sum\limits_{j=1}^n\omega_jx_j\right) 
\right]\,d\eta_1\\
\times \int\limits_{K^{n-1}}\left( 
\max\limits_{1\le j\le n}|\eta_j|\right)^{-n-\alpha }
d\eta_2\ldots d\eta_n.
\end{multline*}

In order to compute the integral over $K^{n-1}$, we write it in 
the form
$$
\int\limits_{K^{n-1}}\left( 
\max\limits_{1\le j\le n}|\eta_j|\right)^{-n-\alpha }
d\eta_2\ldots d\eta_n=I_1+I_2,
$$
$$
I_1=\int\limits_{\max\limits_{2\le j\le n}|\eta_j|<|\eta_1|}
|\eta_1|^{-n-\alpha }d\eta_2\ldots d\eta_n,
$$
$$
I_2=\int\limits_{\max\limits_{2\le j\le n}|\eta_j|\ge |\eta_1|}
\left( \max\limits_{2\le j\le n}|\eta_j|\right)^{-n-\alpha }
d\eta_2\ldots d\eta_n.
$$

It is well known (see, for example, \cite{T68}) that
$$
\int\limits_{\max\limits_{2\le j\le n}|\eta_j|=q^k}d\eta_2\ldots d\eta_n
=q^{(n-1)k}\left( 1-q^{-n+1}\right) .
$$
Suppose that $|\eta_1|=q^\nu$, $\nu \in \mathbb Z$. Then
$$
I_1=|\eta_1|^{-n-\alpha }\sum\limits_{k=-\infty }^{\nu -1}
q^{(n-1)k}\left( 1-q^{-n+1}\right) =|\eta_1|^{-n-\alpha }
q^{(n-1)(\nu -1)}=q^{-(n-1)}|\eta_1|^{-\alpha -1},
$$
$$
I_2=\sum\limits_{k=\nu }^\infty q^{-k(n+\alpha )}q^{(n-1)k}\left( 
1-q^{-n+1}\right) =\left( 1-q^{-n+1}\right) \sum\limits_{k=\nu }^\infty q^{-k(\alpha 
+1)}=\frac{1-q^{-n+1}}{1-q^{-\alpha -1}}|\eta_1|^{-\alpha -1},
$$
so that
$$
\int\limits_{K^{n-1}}\left( 
\max\limits_{1\le j\le n}|\eta_j|\right)^{-n-\alpha }
d\eta_2\ldots d\eta_n=\frac{1-q^{-n-\alpha }}{1-q^{-\alpha -1}}|\eta_1|^{-\alpha 
-1}.
$$

Therefore
\begin{multline*}
\left( D_x^{\alpha ,n}F\right) (t,x)\\
=\frac{1-q^\alpha }{1-q^{-\alpha -1}}\int\limits_K|\eta_1|^{-\alpha 
-1}\left[ f\left( t+\sum\limits_{j=1}^n\omega_jx_j-\eta_1\right)
-f\left( t+\sum\limits_{j=1}^n\omega_jx_j\right) 
\right]\,d\eta_1=\left( D_t^\alpha F\right) (t,x),
\end{multline*}
which means that $F$ satisfies the equation (16). $\qquad 
\blacksquare$

\section{Cauchy Problems}

{\bf 5.1. Applications of the Radon transform.} Let $\varphi \in 
\mathcal D(K^n)$. We will look for a solution $F(t,x)$ of the equation 
(16) satisfying the initial condition
\begin{equation}
F(0,x)=\varphi (x),\quad x\in K^n,
\end{equation}
or the modified initial condition
\begin{equation}
\left( D_t^{n-1}F\right) (0,x)=\varphi (x),\quad x\in K^n.
\end{equation}
Of course, the conditions (17) and (18) coincide if $n=1$.

Let $\check \varphi (\xi ,s)$ be the Radon transform of the 
initial function $\varphi$. Denote
$$
\Gamma (t,x,u)=\check \varphi (u,t+u\cdot x),\quad t\in K,\ x,u\in 
K^n,\|u\|=1.
$$
Let us consider the functions
$$
F_1(t,x)=(1-q^{-1})^{-1}\int\limits_{\|u\|=1} \left( 
D_t^{n-1}\Gamma\right) (t,x,u)\,d^nu,
$$
$$
F_2(t,x)=(1-q^{-1})^{-1}\int\limits_{\|u\|=1}\Gamma 
(t,x,u)\,d^nu.
$$

\medskip
\begin{teo}
The functions $F_1(t,x)$ and $F_2(t,x)$ are radial in $t$, jointly 
locally constant in $(t,x)$, bounded solutions of the Cauchy 
problem (16), (17) and the modified Cauchy problem (16), (18) 
respectively.
\end{teo}

\medskip
{\it Proof}. It follows from the identity (7) that $\check \varphi (\xi 
,r)$ belongs to $\mathcal D(K)$ in $r$ uniformly with respect to 
$\xi \in K^n$, $\|\xi \|=1$ -- there exists a compact set in $K$, 
outside which $\check \varphi (\xi ,\cdot )$ vanishes, for all the 
above $\xi$, and $\check \varphi (\xi ,r+r')=\check \varphi (\xi 
,r)$ if $|r'|\le q^{-l}$ where $l$ does not depend on $\xi$. This 
means that $\Gamma$ is locally constant in $t,x$, uniformly with 
respect to $u\in K^n$ with $\|u\|=1$. In addition, $\Gamma$ and 
$D_t^{n-1}\Gamma$ are bounded, uniformly with respect to $u$. 
These properties make it possible to change the order of 
integration while $D_t^\alpha F_j$ and $D_x^{\alpha ,n}$ are 
computed. Then Proposition 2 shows that $F_1$ and $F_2$ satisfy 
the equation (16). The initial conditions are satisfied due to the 
Radon inversion formula (9). 

\medskip
In order to check that $F_2(t,x)$ is radial in $t$, we notice that
$$
\check \varphi (\omega \xi ,\omega s)=\check \varphi (\xi 
,s),\quad |\omega |=1,
$$
by virtue of the homogeneity property of $\check \varphi$. 
Therefore $\Gamma (\omega t,x,u)=\check \varphi (u,\omega t+u\cdot 
x)=\check \varphi (\omega^{-1}u,t+(\omega^{-1}u)\cdot x)=\Gamma 
(t,x,\omega^{-1}u)$, so that
$$
F_2(\omega t,x)=(1-q^{-1})^{-1}\int\limits_{\|u\|=1}\Gamma 
(t,x,\omega^{-1}u)\,d^nu=F_2(t,x).
$$
Since the operator $D_t^{n-1}$ commutes with the operator 
$f(t)\mapsto f(\omega t)$, $|\omega |=1$, we find also that $F_1$ 
is radial in $t$. $\qquad \blacksquare$

\medskip
Let us study the solution $F_2(t,x)$ of the modified Cauchy 
problem (16), (18) in a little greater detail. Using the 
connection (7) between the Fourier and Radon transforms we get 
that
\begin{equation}
\int\limits_{\|u\|=1}\Gamma (t,x,u)\,d^nu=\int\limits_K\chi 
(-st)\,ds\int\limits_{\|u\|=1}\chi (-s(u\cdot 
x))\widetilde{\varphi }(su)\,d^nu.
\end{equation}

Next,
$$
\int\limits_{\|u\|=1}\chi (-s(u\cdot x))\widetilde{\varphi 
}(su)\,d^nu=\int\limits_{K^n}\varphi (y)\,d^ny\int\limits_{\|u\|=1}\chi 
(s(u\cdot (y-x)))\,d^nu.
$$
By the well-known integration formula (see, for example, 
\cite{T68}),
$$
\int\limits_{\|u\|=1}\chi (s(u\cdot (y-x)))\,d^nu=\begin{cases}
1-q^{-n}, & \text{if $|s|\cdot \|y-x\|\le 1$;}\\
-q^{-n}, & \text{if $|s|\cdot \|y-x\|=q$;}\\
0, & \text{if $|s|\cdot \|y-x\|>q$,}\end{cases}
$$
so that
\begin{equation}
\int\limits_{\|u\|=1}\chi (-s(u\cdot x))\widetilde{\varphi 
}(su)\,d^nu=\left( 1-q^{-n}\right) \int\limits_{\|y-x\| \le 
|s|^{-1}}\varphi (y)\,d^ny-q^{-n}\int\limits_{\|y-x\|=q|s|^{-1}}
\varphi (y)\,d^ny.
\end{equation}

\medskip
\begin{prop}
Suppose that $\varphi (x)=0$ for $\|x\|>q^N$, and $\varphi 
(y)=\varphi (x)$ if $\|y-x\|\le q^{-\nu }$, $\nu ,N\in \mathbb N$. 
Then $F_2(t+t',x)=F_2(t,x)$, if $|t'|\le q^{-\nu}$, and 
$F_2(t,x)=0$ for $|t|>q^{N+1}$.
\end{prop}

\medskip
{\it Proof}. By (19) and (20),
\begin{equation}
F_2(t,x)=(1-q^{-1})^{-1}\int\limits_K\chi (-st)R(s,x)\,ds
\end{equation}
where
\begin{equation}
R(s,x)=\left( 1-q^{-n}\right) \int\limits_{\|y-x\| \le 
|s|^{-1}}\varphi (y)\,d^ny-q^{-n}\int\limits_{\|y-x\|=q|s|^{-1}}
\varphi (y)\,d^ny.
\end{equation}

If $|s|\ge q^{\nu +1}$, then
\begin{multline*}
R(s,x)=\varphi (x)\left\{ \left( 1-q^{-n}\right) \int\limits_{\|y\| \le 
|s|^{-1}}d^ny-q^{-n}\int\limits_{\|y\|=q|s|^{-1}}d^ny\right\} \\
=\varphi (x)|s|^{-n}\left[ \left( 1-q^{-n}\right) -q^{-n}\cdot q^n
\left( 1-q^{-n}\right) \right] =0,
\end{multline*}
so that
$$
F_2(t,x)=(1-q^{-1})^{-1}\int\limits_{|s|\le q^\nu}\chi 
(-st)R(s,x)\,ds,
$$
which implies the required local constancy in $t$.

Let $|t|>q^{N+1}$. Then there exists such an element $s_0\in K$, 
$|s_0|=q^{-N-1}$, that $\chi (s_0t)\ne 1$. If $\| x\|\le q^N$, then 
$\varphi (y)=0$ for $\|y-x\| >q^N$. Therefore for $|s|<q^{-N}$ the
second summand in the right-hand side of (22) equals zero, while 
the domain of integration in the first summand can be fixed as $\{ 
y\in K^n:\ \|y-x\| \le q^N\}$, if $|s|<q^{-N}$. Therefore $R(s,x)$ 
is constant in $s$ on the set $\{ s\in K:\ |s|<q^{-N}\}$, which 
implies the equality $R(s+s_0,x)=R(s,x)$ for all the values of 
$s$. Making in (21) the change of variables $s=s'+s_0$ we come to 
the identity $F_2(t,x)=\chi (s_0t)F_2(t,x)$, which yields the 
required equality $F_2(t,x)=0$. $\qquad \blacksquare$

\medskip
Note that the local constancy of $F_2$ in $t$ may be interpreted 
as a counterpart of the finite domain of dependence for a 
classical wave equation: if the initial function $\varphi$ is such 
that $\varphi (x)=0$ outside some compact set $C\subset K^n$, then 
$F_2(t,x)=0$ for $x\in K^n\setminus C$, at least on some 
neighbourhood of the origin $t=0$. Meanwhile, the fact that 
$F_2(t,x)$ vanishes, as $|t|$ becomes big enough (for a given 
$\|x\|$), resembles the Huygens principle, the existence of the 
trailing edge of a wave. 

\bigskip
{\bf 5.2. A uniqueness theorem.} Here we consider the uniqueness 
problem in the class of generalized solutions, radial in $t$.

Denote by $\Phi'(K,\Phi'(K^n))$ the set of distributions over the 
test function space $\Phi (K)$, with values in $\Phi'(K^n)$.

\medskip
\begin{teo}
Let $F\in \Phi'(K,\Phi'(K^n))$ be a generalized solution of the 
equation (16), that is
$$
\left\langle \left\langle F,D_t^\alpha \varphi_1\right\rangle 
,\varphi_2\right\rangle =\left\langle \left\langle F,\varphi_1
\right\rangle ,D_x^{\alpha ,n}\varphi_2\right\rangle
$$
for any $\varphi_1\in \Phi (K)$, $\varphi_2\in \Phi 
(K^n)$. If $F$ is radial in $t$, then $F\in \mathcal 
D(K,\Phi'(K^n))$. If, in addition, $F(0,x)=0$ or $\left( 
D_t^{n-1}F\right) (0,x)=0$, then $F(t,x)\equiv 0$.
\end{teo}

\medskip
{\it Proof}. Denote by $\widetilde{F}(t,\cdot )$ the Fourier 
transform of $F$ in the variable $x$; as usual, we abuse the 
notation slightly, writing a distribution in the variable $t$ as a 
function of $t$. For any $\psi \in \Psi (K^n)$ we have
$$
D_t^\alpha \langle \widetilde{F}(t,\cdot ),\psi \rangle =\langle 
\| \xi \|^\alpha \widetilde{F}(t,\xi ),\psi (\xi )\rangle .
$$
If $\supp \psi \subset S_N=\left\{ \xi \in K^n:\ \| \xi \| 
=q^N\right\}$, $N\in \mathbb N$, then
$$
D_t^\alpha \langle \widetilde{F}(t,\cdot ),\psi \rangle =q^{\alpha 
N}\langle \widetilde{F}(t,\cdot ),\psi \rangle .
$$

By Proposition 1, the function $\langle \widetilde{F}(t,\cdot ),\psi 
\rangle$ has the form (13) with $t$ substituted for $x$ and some 
$c\in \mathbb C$. If $\psi \in \Psi (K^n)$, then $\psi$ is a sum 
of a finite number of functions supported on spheres $S_N$. 
Taking, in particular, $\psi =\widetilde{\varphi}$, $\varphi \in 
\Phi (K^n)$, we find that $\langle F(t,\cdot ),\varphi \rangle$ 
belongs to $\mathcal D(K)$ in the variable $t$, for any $\varphi \in 
\Phi (K^n)$.

If $F(0,\cdot )=0$, then also $\widetilde{F}(0,\cdot )=0$. If 
$\psi \in \Psi (K^n)$, $\supp \psi \subset S_N$, then, as we have seen, 
$\left\langle \widetilde{F}(t,\cdot ),\psi \right\rangle$ has the 
form (13), and the assumption $\widetilde{F}(0,\cdot )=0$ implies 
the equality $c=0$, whence $\left\langle \widetilde{F}(t,\cdot ),\psi 
\right\rangle =0$, and $\widetilde{F}(t,\cdot )=0$ (since $\psi$ and 
$N$ are arbitrary), and $F(t,\cdot )=0$.

Next, if a function $u(t)$ has a form (13), then its Fourier 
transform has a form (12), and it is easy to find $\left( 
D^{n-1}u\right) (t)$:
$$
\left( D^{n-1}u\right) (t)=\begin{cases}
c(1-q^{-1})q^{Nn}, & \text{if $|t|\le q^{-N}$;}\\
-cq^{Nn-1}, & \text{if $|t|=q^{-N+1}$;}\\
0, & \text{if $|t|>q^{-N+1}$.}\end{cases}
$$
Repeating the above arguments, we find that the equality
$\left( D^{n-1}F\right) (0,x))$ implies $F(t,x)\equiv 0.\qquad 
\blacksquare$

\medskip
It follows from Lemma 1 that bounded locally constant solutions of 
the equation (16) are generalized solutions of the class 
considered in Theorem 2. Therefore the solutions of the Cauchy 
problems constructed in Theorem 1 are unique in the class of 
radial in $t$, bounded locally constant functions. It is natural 
to see such solutions as {\it classical} solutions of the 
non-Archimedean wave equation (16).

\bigskip
{\bf 5.3. Representation of solutions.} Suppose that $\varphi \in 
\Phi (K^n)$. We will look for a solution belonging to 
$\Phi (K)$ and radial in $t$, for each 
$x\in K^n$, and belonging to $\Phi (K^n)$ in $x$, for each $t\in 
K$. In this framework, we may use the Fourier transform, only we 
should not forget to check that the resulting solution indeed 
satisfies the above requirements.

Let us consider the modified Cauchy problem (16), (18). Suppose 
that $n\ge 2$. Performing the Fourier transform in $x$ we come to 
the problem
\begin{equation}
D_t^\alpha \widetilde{F}_2(t,\xi )-\| \xi \|^{\alpha }\widetilde{F}_2(t,\xi 
)=0,
\end{equation}
\begin{equation}
\left( D_t^{n-1}\widetilde{F}_2\right) (0,\xi )=\widetilde{\varphi 
}(\xi ).
\end{equation}
As we have seen,
$$
\widetilde{F}_2(t,\xi )=\begin{cases}
c(\xi )(1-q^{-1})\| \xi \|, & \text{if $|t|\le \|\xi \|^{-1}$;}\\
-c(\xi )q^{-1}\| \xi \|, & \text{if $|t|=q\|\xi \|^{-1}$;}\\
0, & \text{if $|t|>q\|\xi \|^{-1}$,}\end{cases}
$$
where $c(\xi )\in \mathbb C$, $c(0)=0$; note that for $\xi =0$ it 
follows from (23) that $\widetilde{F}_2(t,0)$ is a constant which 
must equal zero by our assumption that $F_2\in \Phi (K)$ in $t$.

Computing $D_t^{n-1}\widetilde{F}_2$ as above (see the proof of 
Theorem 2) we find that
$$
\left( D_t^{n-1}\widetilde{F}_2\right) (t,\xi )=\begin{cases}
c(\xi )(1-q^{-1})\| \xi \|^n, & \text{if $|t|\le \|\xi \|^{-1}$;}\\
-c(\xi )q^{-1}\| \xi \|^n, & \text{if $|t|=q\|\xi \|^{-1}$;}\\
0, & \text{if $|t|>q\|\xi \|^{-1}$,}\end{cases}
$$
We find from the initial condition (24) that $c(\xi 
)=(1-q^{-1})^{-1}\| \xi \|^{-n}\widetilde{\varphi }(\xi )$, and 
come to the expression
\begin{equation}
\widetilde{F}_2(t,\xi )=\| \xi \|^{-n+1}b(t\xi 
)\widetilde{\varphi }(\xi )
\end{equation}
where
$$
b(z)=\begin{cases}
1, & \text{if $\|z\|\le 1$;}\\
-\frac{1}{q-1}, & \text{if $\|z\|=q$;}\\
0, & \text{if $\|z\|>q$.}\end{cases}
$$

Since $\widetilde{\varphi }\in \Psi (K^n)$, it vanishes on a 
neighbourhood of the origin, and it follows from (25) that
$\widetilde{F}_2\in \Psi (K^n)$ in $\xi$, so that $F_2\in \Phi 
(K^n)$ in $x$. In addition, $F_2,\widetilde{F}_2\in \mathcal D(K)$ 
in $t$, uniformly with respect to $x$ (in the sense of support and 
local constancy), which permits to interchange operations in 
different variables. On the other hand, calculating the Fourier 
transforms we obtain from (25) the following representation of the 
solution $F_2(t,x)$:
\begin{equation}
F_2(t,x)=\left( A*B_t*\varphi \right)(x)
\end{equation}
where the convolution is taken with respect to $x$, 
$A(x)=\dfrac{1-q^{-n+1}}{1-q^{-1}}\|x\|^{-1}$, 
$B_t(x)=|t|^{-n}\widetilde{b}(t^{-1}x)$,
$$
\widetilde{b}(\zeta )=\begin{cases}
\frac{q-q^n}{q-1}, & \text{if $\|\zeta \|\le q^{-1}$;}\\
\frac{q}{q-1}, & \text{if $\|\zeta \|=1$;}\\
0, & \text{if $\|\zeta \|>1$.}\end{cases}
$$

The representation (26) makes it possible, for example, to 
investigate the dependence $\varphi \mapsto F_2(t,\cdot )$ with 
respect to the $L_\varkappa$-norms (for a fixed $t\in K$), 
$1<\varkappa <\infty$.

Note that
$$
\|B_t\|_{L_1(K^n)}=|t|^{-n}\int\limits_{K^n}|
\widetilde{\beta}(t^{-1}x)|\,d^nx=\int\limits_{K^n}|
\widetilde{\beta}(x)|\,d^nx,
$$
and the Young inequality, together with the commutativity of 
convolution, gives
$$
\|F_2(t,\cdot )\|_{L_\varkappa}\le C\|A*\varphi \|_{L_\varkappa}
$$
where $C$ does not depend on $t$. Applying a result regarding the 
Riesz potentials from \cite{T68}, we find that for $1<\varkappa 
<\dfrac{n}{n-1}$,
$$
\|F_2(t,\cdot )\|_{L_\lambda }\le C'\|\varphi \|_{L_\varkappa}
$$
where $\lambda =\dfrac{n\varkappa}{n-\varkappa (n-1)}$, and $C'$ 
does not depend on $t$.

For the Cauchy problem (16), (17) (including the case $n=1$), we 
have $F_1(t,x)=\left( B_t*\varphi \right) (x)$, so that
$$
\|F_1(t,\cdot )\|_{L_\varkappa }\le C\|\varphi \|_{L_\varkappa},
$$
for any $\varkappa \in (1,\infty )$, with a constant $C$ 
independent of $t$.

\end{document}